\documentstyle[12pt,thmsa,sw20lart]{article}
\input{tcilatex}
\input tcilatex
\begin{document}

\author{{\small TEODOR\ OPREA}}
\title{{\large ON A RIEMANNIAN\ INVARIANT OF\ CHEN TYPE}}
\date{}
\maketitle

\medskip In [6] we proved Chen's inequality regarded as a problem of
constrained

maximum. In this paper we introduce a Riemannian invariant obtained

from Chen's invariant, replacing the sectional curvature by the Ricci cur-

vature of $k$-order. This invariant can be estimated, in the case of submani-

folds $M$ in space forms $\widetilde{M}(c)$, varying with $c$ and the mean
curvature of

$M$ in $\widetilde{M}(c)$.

\medskip \bigskip

\begin{center}
{\bf 1. INTRODUCTION}
\end{center}

\medskip \medskip \bigskip

We consider a Riemannian manifold $(M,g)$ of dimension $n$, and we fix the
point $x\in M.$ The scalar curvature is defined by 
\[
\tau =\dsum\limits_{1\leq i<j\leq n}R(e_{i},e_{j},e_{i},e_{j}), 
\]
where $R$ is the Riemann curvature tensor of $(M,g)$ and $%
\{e_{1},e_{2},...,e_{n}\}$ is an orthonormal frame in $T_{x}M$.

Let $L$ be a vector subspace of dimension $k\in [2,n]$ in $T_{x}M$. If $X\in
L$ is a unit vector, and $\{e_{1}^{\prime },e_{2}^{\prime
},...,e_{k}^{\prime }\}$ an orthonormal frame in $L$, with $e_{1}^{\prime
}=X $, we shall denote

\[
\text{Ric}_{L}(X)=\dsum\limits_{j=2}^{k}k(e_{1}^{\prime }\wedge
e_{j}^{\prime }), 
\]
where $k(e_{1}^{\prime }\wedge e_{j}^{\prime })$ is the sectional curvature
given by Sp$\{e_{1}^{\prime },e_{j}^{\prime }\}.$

\ 

Using the Ricci curvature of $k$-order at the point $x\in M$

\[
\theta _{k}(x)=\frac{1}{k-1}\min\limits\Sb L,\text{ }\dim L=k  \\ X\in L,%
\text{ }\left\| X\right\| =1  \endSb \text{Ric}_{L}(X), 
\]
we define the invariant

\[
\begin{array}{c}
\delta _{k}(M):M\rightarrow R, \\ 
\delta _{k}(M)=\tau -\theta _{k}.
\end{array}
\]

\ 

For $k=2$ we have $\delta _{k}(M)=\tau -\min (K)=\delta _{M}$, where $K$ is
the sectional curvature and $\delta _{M}$ is the Chen's invariant.

\medskip \bigskip

\begin{center}
{\bf 2. OPTIMIZATIONS\ ON\ RIEMANNIAN\ MANIFOLDS}
\end{center}

~\quad \bigskip

Let $(N,\widetilde{g})$ be a Riemannian manifold, $(M,g)$ a Riemannian
submanifold of $N$, and $f:N\rightarrow R$ a differentiable function$.$ To
these ingredients we attach the optimum problem\newline

\[
(1)\text{ }\min\limits_{x\in M}f(x). 
\]

Let's remember the result obtained in [6].

\bigskip

T{\footnotesize HEOREM} 2.1.{\sl \ If }$x_{0}\in M${\sl \ is a optimal
solution of the problem }$(1)${\sl , then}

{\sl \ }

i){\sl \ }$($grad $f)(x_{0})${\sl \ }$\in T_{x_{0}}^{\perp }M,$

{\sl \ }

ii){\sl \ the bilinear form }

\[
\begin{array}{c}
\alpha {\sl \ :\ }T_{x_{0}}M\times T_{x_{0}}M{\sl \ }\rightarrow {\sl \ }R,
\\ 
\alpha (X,Y)=\text{Hess}_{f}(X,Y)+\widetilde{g\text{ }}(h(X,Y),(\text{grad f}%
)(x_{0}))
\end{array}
\]
{\sl is positive semidefinite, where }$h${\sl \ is the second fundamental
form of the submanifold }$M${\sl \ in }$N.$

\medskip

R{\footnotesize EMARK}. The bilinear form $\alpha $ is nothing else but Hess$%
_{\left. f\right| M}(x_{0}).$

\medskip \bigskip

\begin{center}
{\bf 3. THE INEQUALITY\ SATISFIED\ BY\ THE\ RIEMANNIAN\ INVARIANT\ }$\delta
_{k}(M)$

\bigskip
\end{center}

B.Y. Chen showed in [1] that the Chen's invariant $\delta _{M}$ of a
Riemannian submanifold in a real space form $\widetilde{M}(c)$ satisfies the
inequality

\[
\delta _{M}\leq \frac{n-2}{2}\{\frac{n^{2}}{n-1}\left\| H\right\|
^{2}+(n+1)c\}, 
\]
where $H$\ is the mean curvature vector of submanifold $M$\ in $\widetilde{M}%
(c)$ and $n\geq 3$ is the dimension of $M$. The equality is attained at the
point $x\in M$\ if and only if there is an orthonormal frame $%
\{e_{1},...,e_{n}\}$\ in $T_{x}M$\ and an orthonormal frame $%
\{e_{n+1},...,e_{m}\}$\ in $T_{x}^{\perp }M$\ in which the Weingarten
operators take the following form\newline
\[
A_{n+1}=\left( 
\begin{array}{lllll}
h_{11}^{n+1} & 0 & 0 & . & 0 \\ 
0 & h_{22}^{n+1} & 0 & . & 0 \\ 
0 & 0 & h_{33}^{n+1} & . & 0 \\ 
\text{.} & \text{.} & \text{.} & \text{.} & \text{.} \\ 
0 & 0 & 0 & . & h_{nn}^{n+1}
\end{array}
\right) , 
\]
with $h_{11}^{n+1}+h_{22}^{n+1}=h_{33}^{n+1}=...=h_{nn}^{n+1}$ and\newline
\[
A_{r}=\left( 
\begin{array}{lllll}
h_{11}^{r} & h_{12}^{r} & 0 & . & 0 \\ 
h_{12}^{r} & -h_{11}^{r} & 0 & . & 0 \\ 
0 & 0 & 0 & . & 0 \\ 
. & . & . & . & . \\ 
0 & 0 & 0 & . & 0
\end{array}
\right) ,\forall \text{ }r\in [n+2,m]. 
\]

\medskip

The invariant $\delta _{k}(M)$ satisfies the same inequality. Indeed,
obviously one has min($K$)$\leq \theta _{k}$, which implies $\delta
_{k}(M)\leq \delta _{M}.$ Therefore 
\[
\delta _{k}(M)\leq \frac{n-2}{2}\{\frac{n^{2}}{n-1}\left\| H\right\|
^{2}+(n+1)c\}. 
\]
We give another proof of this inequality for two reasons: to obtain the
equality case and because this proof is useful in order to obtain a stronger
inequality in Lagrangian case.

\bigskip

T{\footnotesize HEOREM }3.1.{\sl \ Consider (}$\widetilde{M}(c),\widetilde{g}%
)${\sl \ a real space form of dimension }$m,\newline
M\subset \widetilde{M}(c)${\sl \ a submanifold of dimension }$n\geq 3,${\sl %
\ and }$k\in [3,n]${\sl . Then}

\[
\delta _{k}(M)\leq \frac{n-2}{2}\{\frac{n^{2}}{n-1}\left\| H\right\|
^{2}+(n+1)c\}, 
\]
{\sl the equality occurring at the point }$x${\sl \ if and only if there is
an orthonormal frame }$\{e_{1},...,e_{n}\}${\sl \ in }$T_{x}M${\sl \ and an
orthonormal frame }$\{e_{n+1},...,e_{m}\}${\sl \ in }$T_{x}^{\perp }M${\sl \
for which the Weingarten operators take the form}

\[
A_{r}=\left( 
\begin{array}{ccccc}
0 & 0 & 0 & . & 0 \\ 
0 & a^{r} & 0 & . & 0 \\ 
0 & 0 & a^{r} & . & 0 \\ 
. & . & . & . & . \\ 
0 & 0 & 0 & . & a^{r}
\end{array}
\right) ,\forall \text{ }r\in [n+1,m].\ 
\]

\medskip

{\it Proof. }Let us consider the point $x\in M$, $\{e_{1},e_{2},...,e_{n}\}$
an orthonormal frame in $T_{x}M$ and $\{e_{n+1},e_{n+2},...,e_{m}\}$ an
orthonormal frame in $T_{x}^{\perp }M.$ \ 

If $L=$Sp$\{e_{1},e_{2},...,e_{k}\}$, then\newline
(1) Ric$_{L}(e_{1})=\dsum\limits_{i=2}^{k}R(e_{1},e_{i},e_{1},e_{i}).$

\ 

From Gauss' equation we obtain the following relations\newline
(2) $\tau =\frac{n(n-1)}2c+\dsum\limits_{r=n+1}^m\dsum\limits_{1\leq i<j\leq
n}(h_{ii}^rh_{jj}^r-(h_{ij}^r)^2)$ and\newline
(3) $(k-1)c=\dsum\limits_{i=2}^kR(e_1,e_i,e_1,e_i)-\dsum\limits_{r=n+1}^m%
\dsum\limits_{i=2}^k(h_{11}^rh_{ii}^r-(h_{1i}^r)^2)$.

\ 

From (1) and (3), it follows\newline
(4) $\frac{\text{Ric}_{L}(e_{1})}{k-1}=c+\frac{1}{k-1}\dsum%
\limits_{r=n+1}^{m}\dsum%
\limits_{i=2}^{k}(h_{11}^{r}h_{ii}^{r}-(h_{1i}^{r})^{2}).$\ 

\ 

From (2) and (4), we obtain\newline
(5) $\tau -\frac{\text{Ric}_{L}(e_{1})}{k-1}=\frac{(n+1)(n-2)}{2}%
c+\dsum\limits_{r=n+1}^{m}\dsum\limits_{1\leq i<j\leq
n}(h_{ii}^{r}h_{jj}^{r}-(h_{ij}^{r})^{2})-$\newline
$-\frac{1}{k-1}\dsum\limits_{r=n+1}^{m}\dsum%
\limits_{i=2}^{k}(h_{11}^{r}h_{ii}^{r}-(h_{1i}^{r})^{2})=$\newline
$=\frac{(n+1)(n-2)}{2}c+\dsum\limits_{r=n+1}^{m}(\dsum\limits_{1\leq i<j\leq
n}h_{ii}^{r}h_{jj}^{r}-\frac{1}{k-1}h_{11}^{r}\dsum%
\limits_{i=2}^{k}h_{ii}^{r})-$\newline
$-\frac{k-2}{k}\dsum\limits_{i=2}^{k}(h_{1i}^{r})^{2}-\dsum%
\limits_{i=k+1}^{n}(h_{1i}^{r})^{2}-\dsum\limits_{2\leq i<j\leq
n}(h_{ij}^{r})^{2}.$

\medskip

As $k\geq 3$, by using (5), one gets\newline
(6) $\tau -\frac{\text{Ric}_{L}(e_{1})}{k-1}\leq \frac{(n+1)(n-2)}{2}%
c+\dsum\limits_{r=n+1}^{m}(\dsum\limits_{1\leq i<j\leq
n}h_{ii}^{r}h_{jj}^{r}-\frac{1}{k-1}h_{11}^{r}\dsum%
\limits_{i=2}^{k}h_{ii}^{r}).$

For $r\in [n+1,m]$, let us consider the quadratic form

\[
\begin{array}{c}
f_r:R^n\rightarrow R, \\ 
f_r(h_{11}^r,h_{22}^r,...,h_{nn}^r)=\dsum\limits_{1\leq i<j\leq
n}h_{ii}^rh_{jj}^r-\frac 1{k-1}h_{11}^r\dsum\limits_{i=2}^kh_{ii}^r
\end{array}
\]
and the constrained extremum problem

\[
\max f_{r} 
\]
\[
\text{subject to }P\text{: }h_{11}^{r}+h_{22}^{r}+...+h_{nn}^{r}=k^{r}\text{,%
} 
\]
where $k^{r}$ is a real constant.

\ 

The partial derivatives of the function $f_{r}$ are\newline
(7) $\frac{\partial f_{r}}{\partial h_{11}^{r}}=\dsum%
\limits_{i=2}^{n}h_{ii}^{r}-\frac{1}{k-1}\dsum\limits_{i=2}^{k}h_{ii}^{r},$%
\newline
(8) $\frac{\partial f_{r}}{\partial h_{jj}^{r}}=\dsum\limits_{i\in \overline{%
1,n}\backslash \{j\}}h_{ii}^{r}-\frac{1}{k-1}h_{11}^{r},$ $j\in [2,k],$%
\newline
(9) $\frac{\partial f_{r}}{\partial h_{ll}^{r}}=\dsum\limits_{i\in \overline{%
1,n}\backslash \{l\}}h_{ii}^{r}$, $l\in [k+1,n].$

\ 

For a optimal solution $(h_{11}^{r},h_{22}^{r},...,h_{nn}^{r})$ of the
problem in question, the vector $($grad) $(f_{1})$ is normal at $P$ that is,
it is colinear with the vector $\left( 1,1,...,1\right) .$

\ 

From (7), (8), (9), it follows that a critical point of the considered
problem has the form\newline
(10) $(h_{11}^{r},h_{22}^{r},...,h_{nn}^{r})=(0,a^{r},a^{r},...,a^{r})$.

As $\dsum\limits_{j=1}^{n}h_{jj}^{r}=k^{r}$, by using (10), we obtain $%
(n-1)a^{r}=k^{r}$, therefore\newline
(11) $a^{r}=\frac{k^{r}}{n-1}$.

Let $p\in P$ be an arbitrary point.

\ 

The 2-form $\alpha $ : $T_{p}P\times T_{p}P\rightarrow R$ has the expression

\[
\alpha (X,Y)=\text{Hess}_{f_{r}}(X,Y)+\left\langle h^{\prime }(X,Y),(\text{%
gradf}_{\text{r}})(p)\right\rangle , 
\]
where $h^{\prime }$ is the second fundamental form of $P$ in $R^{n}$ and $%
\left\langle \text{ },\text{ }\right\rangle $ is the standard inner-product
of $R^{n}$.

\ 

In the standard frame of $R^n$, the Hessian of $f_r$ has the matrix

\[
\text{Hess}_{f_{r}}=\left( 
\begin{array}{cccccccc}
0 & \frac{k-2}{k-1} & \frac{k-2}{k-1} & . & \frac{k-2}{k-1} & 1 & . & 1 \\ 
\frac{k-2}{k-1} & 0 & 1 & . & 1 & 1 & . & 1 \\ 
\frac{k-2}{k-1} & 1 & 0 & . & 1 & 1 & . & 1 \\ 
. & . & . & . & . & . & . & . \\ 
\frac{k-2}{k-1} & 1 & 1 & . & 0 & 1 & . & 1 \\ 
1 & 1 & 1 & . & 1 & 0 & . & 1 \\ 
. & . & . & . & . & . & . & . \\ 
1 & 1 & 1 & . & 1 & 1 & . & 0
\end{array}
\right) \text{.} 
\]

\ 

As $P$ is totally geodesic in $R^{n}$, considering a vector $X$ tangent at $%
p $ to $P$, that is verifying the relation $\sum\limits_{i=1}^{n}X^{i}=0,$
we have

$\alpha (X,X)=-\frac{1}{k-1}[%
(X^{1}+X^{2})^{2}+(X^{1}+X^{3})^{2}+...+(X^{1}+X^{k})^{2}+$\newline
$+(k-2)(X^{2})^{2}+(k-2)(X^{3})^{2}+...+(k-2)(X^{k})^{2}+(k-1)\dsum%
\limits_{i=k+1}^{n}(X^{i})^{2}]\leq 0$. So Hess$_{\left. f\right| P}$ is
negative definite.

Consequently $(0,a^{r},...,a^{r})$, with $a^{r}=\frac{k^{r}}{n-1},$ is a
global maximum point, therefore\newline
(12) $f_{r}\leq \frac{(n-1)(n-2)}{2}(a^{r})^{2}=\frac{(n-2)}{2(n-1)}%
(k^{r})^{2}=\frac{(n-2)}{2(n-1)}n^{2}(H^{r})^{2}$.

\ 

From (6) and (12), it follows that $\tau -\frac{Ric_{L}(e_{1})}{k-1}\leq 
\frac{(n+1)(n-2)}{2}c+$\newline
$+\dsum\limits_{r=n+1}^{m}\frac{(n-2)}{2(n-1)}n^{2}(H^{r})^{2}=\frac{%
(n+1)(n-2)}{2}c+\frac{n^{2}(n-2)}{2(n-1)}\left\| H\right\| ^{2}=\frac{n-2}{2}%
\{\frac{n^{2}}{n-1}\left\| H\right\| ^{2}+$\newline
$+(n+1)c\}$, therefore

\medskip

(13) $\delta _{k}(M)\leq \frac{n-2}{2}\{\frac{n^{2}}{n-1}\left\| H\right\|
^{2}+(n+1)c\}$.

\ 

In (13) we have equality if and only if the same thing occurs in the
inequality (6) and, in addition, (10) occurs. Therefore in (13) we have
equality if and only if there is an orthonormal frame $\{e_{1},...,e_{n}\}$
in $T_{x}M$ and an orthonormal frame $\{e_{n+1},...,e_{m}\}$ in $%
T_{x}^{\perp }M$ for which the Weingarten operators have the following form

\[
A_{r}=\left( 
\begin{array}{ccccc}
0 & 0 & 0 & . & 0 \\ 
0 & a^{r} & 0 & . & 0 \\ 
0 & 0 & a^{r} & . & 0 \\ 
. & . & . & . & . \\ 
0 & 0 & 0 & . & a^{r}
\end{array}
\right) ,\forall \text{ }r\in [n+1,m].\ 
\]

\medskip

\begin{center}
{\bf 4. THE\ LAGRANGIAN\ CASE}
\end{center}

\medskip \bigskip

Let $(\widetilde{M},\widetilde{g},J)$ be a K\"{a}hler manifold of real
dimension $2m.$ A submanifold $M$ of dimension $n$ of $(\widetilde{M},%
\widetilde{g},J)$ is called a totally real submanifold if for any point $x$
in $M$ the relation $J(T_{x}M)\subset T_{x}^{\perp }M$ holds.

If, in addition, $n=m,$ then $M$ is called Lagrangian submanifold. For a
Lagrangian submanifold, the relation $J(T_{x}M)=T_{x}^{\perp }M$ occurs.

A K\"{a}hler manifold with constant holomorphic sectional curvature is
called a complex space form and is denoted by $\widetilde{M}(c)$. The
Riemann curvature tensor $\widetilde{R}$ of $\widetilde{M}(c)$ satisfies the
relation

$\widetilde{R}(X,Y)Z=\frac{c}{4}\{\widetilde{g}(Y,Z)X-\widetilde{g}(X,Z)Y+%
\widetilde{g}(JY,Z)JX-\widetilde{g}(JX,Z)JY+$\newline
$+2\widetilde{g}(X,JY)JZ\}.$

\medskip

{\it Remark}{\bf . }i) If $M$ is a totally real submanifold of real
dimension $n$ in a complex space form $\widetilde{M}(c)$ of real dimension $%
2m,$ then

\medskip 
\[
A_{JY}X=-Jh(X,Y)=A_{JX}Y, 
\]
where $X$ and $Y$ are two arbitrary vector fields on $M.$

ii) Let $m=n$ ($M$ is Lagrangian in $\widetilde{M}(c)$). If we consider the
point $x\in M$, the orthonormal frames $\{e_{1},...,e_{n}\}$ in $T_{x}M$ and 
$\{Je_{1},...,Je_{n}\}$ in $T_{x}^{\perp }M$, then

\medskip 
\[
h_{jk}^{i}=h_{ik}^{j},\forall \text{ }i,j,k\in [1,n], 
\]
where $h_{jk}^{i}$ is the component after $Je_{i}$ of the vector $%
h(e_{j},e_{k}).$

\bigskip

T{\footnotesize HEOREM} 4.1 {\sl Let }$M${\sl \ be a totally real
submanifold of dimension }$n,$ $n\geq 3$ {\sl in complex space form }$%
\widetilde{M}(c)${\sl \ of real dimension }$2m.${\sl \ Then} 
\[
\delta _{k}(M)\leq \frac{n-2}{2}\{\frac{n^{2}}{n-1}\left\| H\right\|
^{2}+(n+1)\frac{c}{4}\}, 
\]
{\sl the equality occurring if and only if there is an orthonormal frame }$%
\{e_{1},...,e_{n}\}${\sl \ in }$T_{x}M${\sl \ and an orthonormal frame }$%
\{e_{n+1},...,e_{2m}\}${\sl \ in }$T_{x}^{\perp }M${\sl \ for which the
Weingarten operators take the form}

\[
A_{r}=\left( 
\begin{array}{ccccc}
0 & 0 & 0 & . & 0 \\ 
0 & a^{r} & 0 & . & 0 \\ 
0 & 0 & a^{r} & . & 0 \\ 
. & . & . & . & . \\ 
0 & 0 & 0 & . & a^{r}
\end{array}
\right) ,\forall \text{ }r\in [n+1,2m].\ 
\]

\medskip

{\it Proof.} Similar with the proof of theorem 3.1.

\medskip \medskip

If $k=n$, and $M$ is a Lagrangian submanifold in the complex space form $%
\widetilde{M}(c)$, the previous result can be improved.

\ \bigskip

T{\footnotesize HEOREM} 4.2{\sl \ Let }$M${\sl \ be a Lagrangian submanifold
in complex space form }$\widetilde{M}(c)${\sl \ of real dimension }$2n${\sl %
, }$n\geq 3.${\sl \ Then }

\[
\delta _{n}(M)\leq \frac{(n+1)(n-2)}{8}c+\frac{(3n-1)(n-2)n^{2}}{2(3n+5)(n-1)%
}\left\| H\right\| ^{2}. 
\]

\ 

{\it Proof}. Let us consider the point $x\in M$, the vector $X\in T_{x}M$
and $\{e_{1},e_{2},...,e_{n}\}$ an orthonormal frame in $T_{x}M$, with $%
e_{1}=X$ $.$ The fact that $M$ is a Lagrangian submanifold imply that $%
\{Je_{1},Je_{2},...,Je_{n}\}$ is an orthonormal frame in $T_{x}^{\perp }M$ .

If $L=T_{x}M$, we shall denote Ric$(X)=$Ric$_{L}(X).$

With an similar argument to those in the previous theorem, we obtain\newline
(1) $\tau -\frac{\text{Ric}(X)}{n-1}=\frac{(n+1)(n-2)}{8}c+\dsum%
\limits_{r=1}^{n}\dsum\limits_{1\leq i<j\leq
n}(h_{ii}^{r}h_{jj}^{r}-(h_{ij}^{r})^{2})-$\newline
$-\frac{1}{n-1}\dsum\limits_{r=1}^{n}\dsum%
\limits_{j=2}^{n}(h_{11}^{r}h_{jj}^{r}-(h_{1j}^{r})^{2})\leq $\newline
$\leq \frac{(n+1)(n-2)}{8}c+\dsum\limits_{r=1}^{n}(\dsum\limits_{1\leq
i<j\leq n}h_{ii}^{r}h_{jj}^{r}-\frac{1}{n-1}\dsum\limits_{r=1}^{n}\dsum%
\limits_{j=2}^{n}h_{11}^{r}h_{jj}^{r})-$\newline
$-\dsum\limits_{1\leq i<j\leq n}(h_{ij}^{i})^{2}-\dsum\limits_{1\leq i<j\leq
n}(h_{ij}^{j})^{2}+\frac{1}{n-1}(\dsum\limits_{j=2}^{n}(h_{1j}^{1})^{2}+%
\dsum\limits_{j=2}^{n}(h_{1j}^{j})^{2}).$

\ 

Using the symmetry in the three indexes of $h_{ij}^{k},$ one gets\newline
(2) $\tau -\frac{\text{Ric}(X)}{n-1}\leq \frac{(n+1)(n-2)}{8}%
c+\dsum\limits_{r=1}^{n}\dsum\limits_{1\leq i<j\leq n}h_{ii}^{r}h_{jj}^{r}-%
\frac{1}{n-1}\dsum\limits_{r=1}^{n}\dsum%
\limits_{j=2}^{n}h_{11}^{r}h_{jj}^{r}-$\newline
$-\dsum\limits_{1\leq i\neq j\leq n}(h_{jj}^{i})^{2}+\frac{1}{n-1}%
(\dsum\limits_{j=2}^{n}(h_{11}^{j})^{2}+\dsum%
\limits_{j=2}^{n}(h_{jj}^{1})^{2}).$

\ 

\medskip

Let us consider the quadratic forms $f_{1},$ $f_{r}:R^{n}\rightarrow R,$ $%
r\in [2,n],$ defined respectively by

$\ $

$f_1(h_{11}^1,h_{22}^1,...,h_{nn}^1)=\dsum\limits_{1\leq i<j\leq
n}h_{ii}^1h_{jj}^1-\frac
1{n-1}\dsum\limits_{j=2}^nh_{11}^1h_{jj}^1-\dsum\limits_{j=2}^n(h_{jj}^1)^2+%
\frac 1{n-1}\dsum\limits_{j=2}^n(h_{jj}^1)^2,$

$\ $

$f_r(h_{11}^r,h_{22}^r,...,h_{nn}^r)=\dsum\limits_{1\leq i<j\leq
n}h_{ii}^rh_{jj}^r-\frac
1{n-1}\dsum\limits_{j=2}^nh_{11}^rh_{jj}^r-\dsum\limits\Sb 1\leq j\leq n  \\ %
j\neq r  \endSb (h_{jj}^r)^2+\frac 1{n-1}(h_{11}^r)^2.$

\ 

We start with the problem 
\[
\max f_{1} 
\]

\[
\text{subject to }P:h_{11}^{1}+h_{22}^{1}+...+h_{nn}^{1}=k^{1}\text{,} 
\]
where $k^{1}$ is a real constant.

\ 

The first two partial derivatives of the quadratic form $f_{1}$ are\newline
(3) $\frac{\partial f_{1}}{\partial h_{11}^{1}}=\dsum\limits_{2\leq j\leq
n}h_{jj}^{1}-\frac{1}{n-1}\dsum\limits_{j=2}^{n}h_{jj}^{1},$\newline
(4) $\frac{\partial f_{1}}{\partial h_{22}^{1}}=\dsum\limits\Sb 1\leq j\leq
n  \\ j\neq 2  \endSb h_{jj}^{1}-\frac{1}{n-1}h_{11}^{1}-2h_{22}^{1}+\frac{2%
}{n-1}h_{22}^{1}.$

\ 

As for a optimal solution $(h_{11}^{1},h_{22}^{1},...,h_{nn}^{1})$ of the
problem in question, the vector grad$(f_{1})$ is colinear with the vector $%
\left( 1,1,...,1\right) ,$ we obtain\newline
(5) $\dsum\limits_{1\leq j\leq n}h_{jj}^{1}-h_{11}^{1}-\frac{1}{n-1}%
\dsum\limits_{j=2}^{n}h_{jj}^{1}=\dsum\limits_{1\leq j\leq
n}h_{jj}^{1}-h_{22}^{1}-\frac{1}{n-1}h_{11}^{1}-2h_{22}^{1}+\frac{2}{n-1}%
h_{22}^{1},$ therefore\newline
(6) $\frac{n-2}{n-1}h_{11}^{1}=\frac{3n-5}{n-1}h_{22}^{1}-\frac{1}{n-1}%
\dsum\limits_{j=2}^{n}h_{jj}^{1}.$

\ 

Similarly we obtain\newline
(7) $\frac{n-2}{n-1}h_{11}^{1}=\frac{3n-5}{n-1}h_{ii}^{1}-\frac{1}{n-1}%
\dsum\limits_{j=2}^{n}h_{jj}^{1}$, $\forall $ $i\in [2,n]$, whence\newline
(8) $h_{22}^{1}=h_{33}^{1}=...=h_{nn}^{1}=a^{1}.$

\ 

The relations (6) and (8) imply\newline
(9) $\frac{n-2}{n-1}h_{11}^{1}=\frac{3n-5}{n-1}a^{1}-a^{1},$ therefore $%
(n-2) $ $h_{11}^{1}=(2n-4)a^{1}$, whence\newline
(10) $h_{11}^{1}=2a^{1}.$

\ 

As $h_{11}^{1}+h_{22}^{1}+h_{33}^{1}+...+h_{nn}^{1}=k^{1}$, by using (8) and
(10), we obtain\newline
(11) $2a^{1}+(n-1)a^{1}=k^{1}$, therefore\newline
(12) $a^{1}=\frac{k^{1}}{n+1}$.

\medskip \medskip

As $f_{1\text{ }}$is obtained from the function studied in theorem 3.1 by
subtracting some square terms, $f_{1}\left| P\right. $ will have the Hessian
negative definite. Consequently the point $%
(h_{11}^{1},h_{22}^{1},...,h_{nn}^{1})$ given by the relations (8), (10) and
(12) is a maximum point, and hence\newline
(13) $f_{1}\leq 2a^{1}(n-1)a^{1}+C_{n-1}^{2}(a^{1})^{2}-\frac{1}{n-1}%
2a^{1}(n-1)a^{1}-(n-1)(a^{1})^{2}+$\newline
$+\frac{1}{n-1}(n-1)(a^{1})^{2}=\frac{(a^{1})^{2}}{2}(n^{2}-n-2)=\frac{%
(a^{1})^{2}}{2}(n+1)(n-2).$

\medskip

From (12) and (13), one gets\newline
(14) $f_{1}\leq \frac{(k^{1})^{2}}{2(n+1)}(n-2)=\frac{(n-2)n^{2}}{2(n+1)}%
(H^{1})^{2}.$

\ 

Further on, we shall consider the problem 
\[
\max f_{2} 
\]
\[
\text{subject to }P:h_{11}^{2}+h_{22}^{2}+...+h_{nn}^{2}=k^{2}\text{,} 
\]
where $k^{2}$ is a real constant.

\ 

The first three partial derivatives of the quadratic form $f_{2}$ are\newline
(15) $\frac{\partial f_{2}}{\partial h_{11}^{2}}=\dsum%
\limits_{j=2}^{n}h_{jj}^{2}-\frac{1}{n-1}\dsum%
\limits_{j=2}^{n}h_{jj}^{2}-2h_{11}^{2}+\frac{2}{n-1}h_{11}^{2},$\newline
(16) $\frac{\partial f_{2}}{\partial h_{22}^{2}}=\dsum\limits\Sb 1\leq j\leq
n  \\ j\neq 2  \endSb h_{jj}^{2}-\frac{1}{n-1}h_{11}^{2},$\newline
(17) $\frac{\partial f_{2}}{\partial h_{33}^{2}}=\dsum\limits\Sb 1\leq j\leq
n  \\ j\neq 3  \endSb h_{jj}^{2}-\frac{1}{n-1}h_{11}^{2}-2h_{33}^{2}.$

\ 

For a solution $(h_{11}^{2},h_{22}^{2},...,h_{nn}^{2})$ of the problem in
question the vector grad$(f_{2})$ is colinear with $(1,1,...,1)$.

Consequently $\dsum\limits_{j=1}^{n}h_{jj}^{2}-h_{22}^{2}-\frac{1}{n-1}%
h_{11}^{2}=$ $\dsum\limits_{j=1}^{n}h_{jj}^{2}-h_{33}^{2}-\frac{1}{n-1}%
h_{11}^{2}-2h_{33}^{2}$, therefore\newline
(18) $h_{22}^{2}=3h_{33}^{2}.$

\medskip

Similarly we obtain\newline
(19) $h_{22}^{2}=3h_{jj}^{2}=3a^{2}$, $\forall $ $j\in [3,n]$.

\medskip

From (15), (16) and (19) we obtain\newline
(20) $3h_{11}^{2}-\frac{3}{n-1}h_{11}^{2}=3a^{2}-\frac{1}{n-1}%
(3a^{2}+(n-2)a^{2})$, hence\newline
(21) $h_{11}^{2}=\frac{2a^{2}}{3}.$

\ 

We shall denote $a^{2}=3b^{2}$. The relations (19) and (21) becomes\newline
(22) $h_{11}^{2}=2b^{2},$\newline
(23) $h_{22}^{2}=9b^{2},$\newline
(24) $h_{33}^{2}=...=h_{nn}^{2}=3b^{2}.$

\ 

As $h_{11}^{2}+h_{22}^{2}+...+h_{nn}^{2}=k^{2},$ we obtain $%
2b^{2}+9b^{2}+(n-2)3b^{2}=k^{2}$, therefore\newline
(25) $b^{2}=\frac{k^{2}}{3n+5}.$

\ 

With an similar argument to those in the previous problem we obtain that the
point $(h_{11}^{2},h_{22}^{2},...,h_{nn}^{2})$ given by the relations (22),
(23), (24) and (25) is a maximum point. Hence\newline
(26) $f_{2\text{ }}\leq
2b^{2}(9b^{2}+(n-2)3b^{2})+9b^{2}(n-2)3b^{2}+C_{n-2}^{2}(3b^{2})^{2}-$%
\newline
$-\frac{1}{n-1}2b^{2}(9b^{2}+(n-2)3b^{2})-(2b^{2})^{2}-(n-2)(3b^{2})^{2}+%
\frac{1}{n-1}(2b^{2})^{2}=$\newline
$=\frac{(b^{2})^{2}}{2(n-1)}(9n^{3}-6n^{2}-29n+10)=\frac{(b^{2})^{2}}{2(n-1)}%
(3n+5)(3n-1)(n-2).$

\ 

From (25) and (26) we obtain $f_{2\text{ }}\leq \frac{(k^{2})^{2}(3n-1)(n-2)%
}{2(3n+5)(n-1)}=\frac{(3n-1)(n-2)n^{2}}{2(3n+5)(n-1)}(H^{2})^{2}.$

\ 

Similarly one gets\newline
(27) $f_{r}\leq \frac{(3n-1)(n-2)n^{2}}{2(3n+5)(n-1)}(H^{r})^{2}$, $\forall $
$r\in [2,n]$.

\medskip

As $\frac{n-2}{n+1}\leq \frac{(3n-1)(n-2)}{(3n+5)(n-1)}$ , $\forall $ $n\geq
3,$ from (14) and (27), we obtain\newline
(28) $f_{r\text{ }}\leq \frac{(3n-1)(n-2)n^{2}}{2(3n+5)(n-1)}(H^{r})^{2}$ , $%
\forall $ $r\in [1,n].$

\ 

From (2) and (28), one gets\newline
(29) $\tau -\frac{\text{Ric}(X)}{n-1}\leq \frac{(n+1)(n-2)}{8}%
c+\dsum\limits_{r=1}^{n}\frac{(3n-1)(n-2)n^{2}}{2(3n+5)(n-1)}(H^{r})^{2}=$%
\newline
$=\frac{(n+1)(n-2)}{8}c+\frac{(3n-1)(n-2)n^{2}}{2(3n+5)(n-1)}\left\|
H\right\| ^{2}$, therefore\newline
(30) $\delta _{n}(M)\leq \frac{(n+1)(n-2)}{8}c+\frac{(3n-1)(n-2)n^{2}}{%
2(3n+5)(n-1)}\left\| H\right\| ^{2}.$

\end{document}